\documentclass[12pt]{article}
\usepackage{latexsym}
\usepackage{amsmath,amsfonts,amssymb,enumerate,multicol,tikz}
\usepackage{epsfig}
\usepackage{pstricks}
\usepackage{hyperref}
\usepackage[utf8]{inputenc}

\usepackage{mathrsfs}

\def\square{\Box}
\begin{document}
\title{Around Logical Perfection} 
\author{J. A. Cruz Morales \\
A. Villaveces \\
Universidad Nacional de Colombia\\
B. Zilber \\ 
University of Oxford}

\date{}

\maketitle

\begin{abstract}In this article we present and describe a notion of "logical perfection". We extract the
    notion of "perfection" from the contemporary logical concept of categoricity. Categoricity (in
    power) has become in the past half century a main driver of ideas in model theory, both
    mathematically (stability theory may be regarded as a way of approximating categoricity) and
    philosophically. In the past two decades, categoricity notions have started to overlap with more
    classical notions of robustness and smoothness. These have been crucial in various parts of
    mathematics since the nineteenth century. 

We postulate and present the category of logical perfection. We draw on various notions of
perfection from mathematics of the 19th and 20th centuries and then trace the relation to the
concept of categoricity in power as a logical notion of what a "mathematically perfect" structure
is.
\end{abstract}

This essay is an attempt to present the idea of \textbf{logical perfection} to
a philosophical audience. This expression is often used informally in
mathematical practice and sometimes also in more formal discussion around
mathematics. This happens sometimes in the form of an aesthetic criterion, and
is one of the strongest drivers of mathematical activity and one of the main
tests for its relevance.
Since the advent of the discipline of mathematical logic it has become possible
to  investigate a potentially adequate formal notion by mathematical means.
Roughly, before giving more detail, \emph{a mathematical object of a certain
``size''  is logically perfect if in a certain formal language it allows a
``concise'' description fully determining the object.}\footnote{The essay has
  arisen from many conversations and collaboration the authors have had
  during the past few years and it is originally based on two talks the first
  author gave, one in Paris and the other one in Bogotá, about the work of the
  third author. The second author then brought some additional perspective.
  The first author wants to mention in particular the talk given at Bogotá
  during the workshop \emph{Mapping traces: Representation from Categoricity to
  Definability} organised by the second author
  and María Clara Cortés at Universidad Nacional de Colombia in 2014. This was
  very helpful since the audience consisting of philosophers, mathematicians and
  artists made the idea of writing about logical perfection for a general
audience possible.}
   
This notion, in particular, is central in the third author's paper~\cite{zil1}
and has been implicitly present, mainly as a motivating factor, in a number of
other research papers in various branches of mathematics.

%We dare to present this concept to the big public since we think that aside
%from its mathematical relevancy the discussion about the idea of logical
%perfection along the lines we have undertaken in this text is also of interest
%for philosophers and in general for anyone interested in knowing about what
%mathematicians are working on and how this can influence their own knowledge
%and practices, hopefully artists belong to this group too. Therefore, we want
%to show to the layman a good idea about some of the motivations, objects and
%structures present in our work and reflections as mathematicians.

%On the other hand, we want to make special emphasis on the idea of perfection since its search is an important feature of the mathematical activity as the {F}rench philosopher 
%Albert Lautman has pointed out\footnote{We thank Fernando Zalamea for calling our attention to this aspect of Lautman's philosophy. 
%In some sense, we can say that the efforts to clarify the idea of logical perfection follow Lautman's 
%spirit.}. In the next section we will back to this point with more detail. For now, it is worth noting that the search of perfection also leads { many} other human enterprises. \\ 

Of course, writing for a wider audience means we skip many subtle
mathematical details and avoid as much as possible using technical terms. 

We draw on various analogies to show that logical perfection has strong
versions outside of mathematics. Moreover, we will argue that logically
perfect structures can be used for the study of the physical world, making the
idea
relevant not just in mathematics but in the realm of physics too. The question
whether logical perfection manifests itself in other areas of the human
activity (such as art) is left open here; we may only hope it will raise the interest of some of
our readers.

%This text is also a \emph{manifesto.} 
%A manifesto in which we propose a debate  (not restricted to the mathematical community) 
%on the (delicate) notion of logical perfection, and also a manifesto of what
%we consider is one of the future (and fruitful) lines of research in our
%area. \\
%
Finally, we thank the referee for many insightful comments that led to serious
improvements and clarifications of this paper.

\section{Why logical perfection?}\label{whyperf}

The interest in looking for some kind of \emph{perfection} in mathematical structures
is not new.
In the history of their discipline,  mathematicians have been driven to
think about this kind of perfection, albeit for different sorts of motivations,
and have tried to capture this idea by means of
mathematical tools.
Let us mention a few of these attempts in the work of Galois, Riemann
and Grothendieck and build up a first collection of examples for our discussion
of logical perfection.

Galois made a bold switch from the classical perspective of looking directly for
solutions to algebraic equations to a study of the \emph{symmetry} of possible
solutions: a move toward a \emph{completion} of the set of \emph{possible
solutions} by means of the study of the group of symmetries of all the
solutions that could exist, and by filtering out the interaction between
enlarging the field where these solutions could appear and the group of such
symmetries. The resulting theory (aptly named much later ``Galois Theory'')
goes way beyond the initial quest for solutions to algebraic equations and
changed the ground to an idea (the symmetry of possible solutions in extensions
of fields, and the duality between the groups of these symmetries and the
field extensions) that is still central after two centuries, and whose scope
goes much beyond the wildest dreams Galois could have had. This is a first form
of ``perfection'' for us: \textbf{completeness of possible solutions,} and \textbf{register
of emerging symmetry.} More importantly even, underlying these two aspects of perfection arising
in Galois' work, there is a kind of \emph{uniqueness}, only reached once all possible solutions
(and their symmetries) are considered. This wholeness, this uniqueness, albeit implicit in the work
of Galois, is a component of the main tenets of the notion of \emph{logical} perfection we propose.

Our second example is Riemann's work on the foundations of geometry. In a move
parallel to Galois', he went beyond understanding geometry in terms of global
axioms and laid the ground for a ``local'' approach, driven by a ``metric'' (a
way of measuring distances) that could change, twist, curve \emph{itself}.
Instead of pla\-cing objects such as curves, planes, surfaces inside a ``global
space'' (as had been done for aeons in mathematics), Riemann put the
\emph{twisting} itself, so to speak, at centerstage: instead of placing
twisted, curved objects ``in'' a space, the space itself became the twisting.
Here, the notion of logical perfection is of a different kind from what we had
in the Galois example. There, global symmetry (and the connection between
symmetry and extensions where solutions live) was the expression of that
perfection. Here the perfection is rather the new \emph{flexibility} Riemann's
construction offers us, when compared with earlier incarnations of ``space'',
of geometry\footnote{Suffice it to say that half a century later, Einstein would base
his General Relativity Theory on Riemann's work: the mathematical content of
Einstein's theory is essentially present in Riemann's approach to geometry.
Here, the ``perfection'' aspect has more the flavor of a way to construct
\emph{many} possible geometries, one for each ``Riemann metric'' ---one for
each way to ``twist'' space, so to speak--- and a global treatment of all of
these geometries (and moreover, mathematical ways of classifying and comparing
them).}.

Our third example, much more contemporary and of a different kind, is
Grothendieck's new foundations of algebraic geometry. Very roughly, the concept
of a general notion of ``space'' is again at stake. But here Grothendieck
essentially first ``disassembles'' the surfaces or curves (called more
generally ``varieties'' by mathematicians) by putting all the weight of the
analysis into one single aspect (localisation) of the space and then finding a
system for placing these localisations in a coherent way. By doing this,
Grothendieck creates a version of ``space'' (called \emph{affine scheme}) that
embodies two movements: first, the localisation (and the possibility of
treating only \emph{one aspect} of the space) and second, the coherence. This
highlights yet a different aspect of logical ``perfection'': the possibility of
regarding space as a coherent way of pasting localised versions of itself.

What is new in the approach presented here is that we claim the existence of a
relevant rigorous mathematical concept which allowed an amazingly deep theory,
and has led to a new understanding of a number of specific  structures central
to modern mathematics. This rigorous concept is defined within the discipline
called {\bf model theory}, a subdiscipline of mathematical logic, which deals
with formal languages  and their semantics. One can confidently claim that the
central concept of present-day model theory is that of {\bf stability} of
formal theories and one key notion of stability theory (from which it started
in the 1960s) is that of {\bf uncountably categorical} theories. Through the
efforts of many people, and most prominently by contributions of S.~Shelah (see
\cite{Sh}), we now have a rather comprehensive {\bf classification theory}
which establishes an effective hierarchy in the ``universe'' of mathematical
structures (or their theories)\footnote{An interesting interactive
  visualisation of ``a map of the universe'' can be seen online at 
http://www.forkinganddividing.com .   
Unfortunately, the graphics present it as a flat landscape although there is a
natural feel that the ``more stable'' structures should be at  higher levels of
the landscape.}. The hierarchy is effectively based on the complexity of the
system of invariants  which ultimately describe a given structure, a model of a
formal theory\footnote{Shelah also uses the criteria of whether the given
  first-order theory has a {\em structure theorem}, that is if the isomorphism
types of models of the theory can be classified in terms of a simple
combinatorial structure.}.
The highest level of the hierarchy corresponds to the simplest system of
invariants. This corresponds, in some sense, to a highest level of
perfection.

\medskip

The previous emphasis on the stability hierarchy, and in particular the region
near its ``top'' (uncountably categorical theories) \emph{describes} a
mathematically rigorous (and completely abstract) approach to a notion
relevant to a working definition of logical perfection. We still have to
address the issue of how adequate and useful this notion is, which dividing
lines it draws and which important mathematical structures satisfy the
criteria.

An interesting observation from the cumulated experience with the study of the
stability hierarchy would establish (very roughly) that \emph{the higher a
  structure is in the hierarchy, the closer it is to a ``centre of
mathematical universe''.} We may take this center to be algebraic geometry in
the broadest sense\footnote{Definining exactly what algebraic geometry ``in the
broadest sense'' means is not immediate, but we may take Grothendieck's ideas
as the main guide.}. In some (limited) sense we may define the most general
form of geometry to be the structures populating the top levels of stability
hierarchy\footnote{Working on this presumption one arrives at  a meaningful notion of
non-classical geometric spaces (see \cite{cruzzil}, \cite{zil4}, \cite{zil} and
the discussion in section~\ref{s4} ) which in a more conventional mathematical
setting are treated via the formalism of {\em non-commutative (or quantum)
geometry}. The latter approach  is essentially a syntactic  algebraic analysis
avoiding  geometric semantics.}.

\section{Logical perfection and the issue of uniqueness}\label{sectionUniq}

In the previous section we posited one reason why we may consider
\emph{categoricity} (in uncountable cardinals) as a center of classification
theory: the observation that many ``central mathematical structures'' (those
from algebraic geometry or those corresponding to linear phenomena) seem to
hover close to that region\footnote{There are important exceptions to this
reason. The first one is obvious: real numbers are far from being categorical
yet are also clearly central mathematical structures in many senses. However,
aside from the infinite order that is the reason for their
non-categoricity, the exhibit a rather simple structure of \emph{definable
sets}: each one of them is really a finite union of intervals. This notion,
called \textbf{o-minimality,} provides reasons to place them in a region where
some of the good properties of uncountably categorical structures still work,
albeit in a different way. The role of interactions between complex analysis
and real analysis is mimicked by this correspondence. The second exception is
subtler: classification theory provides many other regions that, while not
corresponding to the ``supremely perfect'' uncountably categorical region, they
exhibit very strong regularity and smoothness properties.}.

The notion of categoricity concretises the meaning of {\em uniqueness}. One
says that a collection of statements in a formal language {\em (set of axioms)}
is \emph{categorical} if it has just one model, up to isomorphism. This
expression ``up to isomorphism'' means that we do not want to distinguish two
structures if they differ only by the way their elements are presented. 

The choice of the formal language is very essential. Usually it is meant to be
a first-order language,
that is one which allows only finite length formulas and quantifiers ``for
all'' and ``there exists'' which refer to elements of the structure in question
(but not to relations or functions). However, as the research in the last three
decades has shown, much of what will be said below about categoricity in the
first-order context holds in a more general setting.

The notion of categoricity has existed for as long as logic has been
formalised. But in the context of first order languages one realises very
quickly, from basic facts of the theory, that the above {\em absolute}
categoricity can only hold for descriptions of finite structures. For  infinite
structures $\mathbf{M}$ it is possible to have uniqueness in some cases  if we
add to the first order description the (non-first-order) statement fixing the
cardinality $\kappa$ of the structure  $\mathbf{M}$. This relative categoricity
is called {\bf categoricity in cardinality (in power) $\kappa$} or
$\kappa$-categoricity.
 \\

One has to distinguish two  types of cardinalities in the context of
categoricity, namely, uncountable (large) and countable (the minimal infinite)
categoricity. We are interested in uncountable categorically describable
structures
which entails that the structure is much bigger than the size of its description.
 A remarkable fact was proved by Michael Morley in 1964, namely, that
 categoricity in one uncountable cardinality implies the categoricity in all
 uncountable cardinalities: 
 the actual value of the uncountable cardinal  is irrelevant\footnote{This is in sharp contrast
 with countable categoricity. Countably categorical structures might also in some sense be
 candidates to a kind of perfection, probably - but all the geometric features of
 uncountably categorical structures are lost in that case. This dependence on the cardinality
 might be regarded as non-logical in some sense, but the case of uncountable categoricity has
 strong logical properties as well as strong geometric properties.}. \\

The study of this kind of structures has been in the focus  of research in model theory for at least 60 last years. The amazing 
conclusion derived from the research is that among the huge diversity of
mathematical structures there are very few which satisfy the 
(slightly narrower) definition of categoricity, and those can be classified.
These certainly seem to corresponding to an ideal of {\em logical perfection},
in the following sense: categorical structures $M$ determine a first order
theory $Th(M)$ (the set of all sentences that are true in $M$) and then comes
the reason why we call them ``logically perfect'': \textbf{all other
structures that satisfy the theory $Th(M)$ and are of the same cardinality as
$M$ are isomorphic to $M$.} In other words, uncountably categorical structures
are inextricably linked to their logical description; the description
$T=Th(M)$
completely determines the structure $M$ (with the usual caveat of ``up to
isomorphism'' and because of limitations in the expressive power of first order
logic\footnote{Namely, the Löwenheim-Skolem theorem.} provided also one
considers only structure of the same cardinality as $M$).
\\

It is not that surprising that a remarkable example of such theory is the
theory of the field of complex numbers $\mathbb{C}$ in the language based on
algebraic operations $+$ and $\times$. Note that this is the language where  
algebraic geometry is naturally done\footnote{In algebraic geometry classical objects
are solution sets of algebraic expressions, that is, polynomials written with $+$ and $\cdot$}  but we can not, e.g. distinguish the real part of a
complex number, so we can not speak about the real numbers when working over
$\mathbb{C}.$ Recall that  the theory of the field $\mathbb{R}$ of real numbers is
not categorical\footnote{And is not even stable!}. 
    
Complex numbers are present everywhere in
mathematics as are the reals. However, there is a significant difference in the theories and in fact complex geometry and the geometry of real manifolds are two different specialisations within mathematics. Classification theory detects the difference and following the above logic in effect claims a certain ``priority'' of complex geometry.  

Of course, for a mathematician the choice of an area of research is a personal
matter and is usually made on either historic or aesthetic grounds.  Both
complex and real geometry are equally respected fields of mathematical research
although from our point of view the first is fundamental while the second is
auxiliary.  We stress again the fact that it was Bernhard Riemann 
who first understood how real and complex
geometries interact  with one another and how the study to the latter
introduces a whole new range of powerful methods of algebraic geometry into
the field.

Different criteria work in the studies of real world.  Here the wrong choice of
mathematical setting can have adverse effect on the understanding of reality.
The mathematical model of Newtonian physics was based on real analytic
geometry. This tradition continued into the new physics with the model enriched
by more and more uses of complex numbers, seen rather as convenient auxiliary
tools. One of the first who pointed to the importance of reversing this
perspective was Roger Penrose in his 1978 address at the International Congress
of Mathematicians under the title ``The complex geometry of the natural
world'',~\cite{RP}. In more recent decades, with the arrival of string theory,
the priority or at least the centrality of complex geometry is  undeniable.

To summarise the \emph{logicality} of our notion of perfection: we started with
various notions of perfection as we did in Section~\ref{whyperf}, coming from
differing examples in the history of mathematics but then in this section we
narrowed our focus to the notion of uniqueness and its logical expression,
(uncountable)
\textbf{categoricity}. Then we remarked that a whole classification theory that
encompasses \emph{all first order theories}\footnote{and in more recent decades
much more} on the one hand grew up out of the attempts to prove the Morley
theorem and its generalisations and on the other hand ended up providing ways
of callibrating \emph{exactly how far} from categoricity one is, in terms of
smoothness/regularity properties that slowly vanish as we go further and
further away from categoricity. It is in this very sense that categoricity has
been playing the role of a logical form of perfection. A posteriori we realise
that a major part (although not all) of \emph{central} mathematics actually
happens to be one of the theories that are uncountably categorical.

\section{Logically perfect structures: the role of geometry} \label{s4}

Perhaps the most remarkable feature of model-theoretic classification theory is
that it exposes a geometric nature of some ``perfect'' structures. The
geometric features of those structures arise from their logical definition,
albeit in a highly non-trivial and initially unforeseen way. These were
discovered in the course of proving the original ground-breaking categoricity
theorem of Michael Morley (see previous section) as the key technical instruments of the
proof: Morley rank, homogeneity and, added in later versions of the proof,
\textbf{dimension} (Baldwin and Lachlan), and associated combinatorial
geometries (Marsh, Zilber). It took a while to realise the \emph{geometric
character} of the technical definitions and to develop a new geometric
intuition around the notions. In particular, Morley rank is a very good
analogue of dimension in algebraic and analytic geometry and thus we can think
of ``curves'', ``surfaces'' and so on in the very general context of
categorical and even stable theories.
This stage of the theory is summarised in the monograph \cite{Pi} by A.~Pillay.

In the 1980s the third author formulated a
 Trichotomy Conjecture (see~\cite{zil2}) which, based on the above intuition,
 suggested that any  uncountably categorical structure is ``reducible'' to
 either an object of algebraic geometry, or linear algebra, or to a simple
 combinatorial structure. Although in many special classes the conjecture has
 been confirmed, the general case was refuted by Ehud Hrushovski who found
 remarkable counter-examples opening fascinating new perspectives on the nature
 of model theory (its interactions with geometry) and its links with the
 analytic world.

Around the same time, a way to fix the Trichotomy conjecture was found. This
required narrowing the class of structures subject to the conjecture ---in some
sense, this amounted to refining the notion of logical perfection. This was
done by being more careful in choosing the logic in question.
Namely, our logical language must be
able to  distinguish {\em positively formulated} statements from their
negations. The axioms of a good (perfect) theory must be  ``equational'' just
like laws of  physics and objects of geometry are given by equations
(and never by negating equations). And this is already the principle on which
algebraic geometry is built on! It studies curves, surfaces, shapes given as
solution sets for systems of algebraic equations. Algebraic geometry treats
such sets as {\em closed in the Zariski topology}. The corresponding
generalisation of this notion in the context of categorical and stable
structures led to the notion of a {\bf Zariski structure} (or  {\em Zariski
geometry})\label{zariskipage} introduced by Hrushovski and the third author.

This improvement in the notion led to a desired Classification Theorem
(Hrushovski, Zilber 1993, see \cite{zil0}):

 The class of Zariski geometries  satisfies the Trichotomy Principle and
 therefore Zariski geometries are reducible to\footnote{Here ``reducible to''
   can be taken in a first reading as a technical
nuisance not requiring much explanation. The typical example of Zariski
geometry is a (complex) algebraic variety (glued from affine charts) with
possibly a vector bundle over it,  a description of which can require quite a
lot of technical detail. Such a description eventually reduces to the structure
of the complex field itself. However, the constructions described by the
theorem can go beyond the technicalities of this example, so beyond algebraic
and complex geometry. Ten years after the classification theorem, a closer
analysis of what ``reducible to'' could  mean led to the discovery that a huge
source of new Zariski structures is non-commutative (or quantum) algebraic
geometry, see \cite{zil4}.
} classical structures such as the
 field of complex numbers and vector spaces. \\

It is hard to describe what exactly the subject of geometry as practised by
mathematicians is, but non-commutative geometry is a much bigger mystery. It is
best identified as the study of algebraic structures,  non-commutative {\em
coordinate rings}, that supposedly correspond to hypothetical geometric spaces
which are not necessarily visualisable. Historically, these were physicists
who, starting from the famous ``magic paper'' of Heisenberg of
1925~\cite{HMagic}, have given
up to the attempts to describe the physics of micro-world in classical terms
and instead used a purely formal algebraic calculus (algebraic quantum
mechanics) to successfully explain the behaviour of elementary particles. 
One can say that the physics of the micro-world lives in an unusual, previously
unknown,  geometric space which requires a non-commutative algebra to describe.
Paralleling   this the very centre of the logical universe is occupied by
structures which mathematically stem from the same source.\\ \\

%In effect, this refinement of logic is a further  sharpening of our geometrical intuition and deeper understanding of the role of geometry in our perception of reality, including mathematical ``reality''. 
 The fusion of geometry with other branches of 
mathematics, for instance, number theory and representation theory, was one of
the biggest programs in the mathematics of the 20th century\footnote{The 
figure of Grothendieck was essential in formulating and developing this program
in the broadest generality.}. We would like to believe that the fusion of logic
(model theory)  with other branches of mathematics  
is one of the biggest and ambitious programs of the mathematical research for
the
21st century. In particular, the ``new geometry'' arising from model 
theoretical considerations has the potential to become an important area of
research in mathematics and beyond. 
And the study of logically perfect structures gives a crucial insight.  \\

Summarizing, the search of logically perfect structures leads to consider
geometric/topological ingredients in logic which has as a consequence
that a refinement of the idea of logical perfection is obtained. During this
process the idea of Zariski structures arises from purely logical 
considerations but with a geometrical flavor and motivation. So far, our
discussion has not left the realm of mathematics but as  
our previous discussions (and the title of the essay suggests) we want to go
beyond mathematics, entering the ``real world''. A question arises: Are  
logically perfect structures helpful for understanding the ``real world''? We
answer this question in the positive and now provide some insights for
exploring 
that possibility.

\section{Logical perfection and physics}

We now focus on a different kind of problem: programs for new foundations of
quantum gravity, and the issue of tackling an appropriate notion of
geometric space for physics. On the face of it, this problem would seem quite
remote from our notion of logical perfection. There is however a deep link, as
we will describe.

Let us quote again Roger Penrose (his ICM address~\cite{RP}):

\begin{quotation}
  \emph{
``Even at the most elementary level, there are still severe conceptual problems
in providing a satisfactory interpretation of quantum mechanical observations
in a way compatible with the tenets of special relativity. And quantum field
theory, which represents the fully special-relativistic version of quantum
theory, though it has had some very remarkable and significant successes,
remains beset with inconsistencies and divergent integrals whose illeffects
have been only partially circumvented. Moreover, the present status of the
unification of general relativity with quantum mechanics remains merely a
collection of hopes, ingenious ideas and massive but inconclusive calculations.
\\
In view of this situation it is perhaps not unreasonable to search for a
different viewpoint concerning the role of geometry in basic physics. Broadly
speaking, "geometry", after all, means any branch of mathematics in which
pictorial representations provide powerful aids to one's mathematical
intuition. It is by no means necessary that these "pictures" should refer just
to a spatio-temporal ordering of physical events in the familiar way...''
}
\end{quotation}

\medskip
Penrose continues to discuss structures of complex geometry as new
geometric tools in quantum physics. However, today this seems to be far from
enough. From a similar reasoning the physicist C.~Isham and the philosopher
of physics J.~Butterfield reached a bold program for building a new
foundation of quantum gravity physics, based on Grothendieck toposes as
the most general form of geometric space (see~\cite{IB}).

Naturally, Isham-Butterfield is not the only program to tackle the problem (see
e,g,  the non-commutative geometry approach \cite{CM} by A.~Connes and
M.~Marcolli, which however does not reveal a geometric space as such)
but it seems to be the most ambitious and general\footnote{Maybe too general as
  to the best of our knowledge there is no interesting calculation produced out
  of it.}.

A project, which may be seen as similar in spirit is suggested and started in
\cite{zil} and in shorter form in  \cite{cruzzil}. Like other such programs the
key is the respective notion of the geometric space for physics. Our suggestion
is based on the philosophy of logical perfection; after all it is reasonable to
expect that the geometric structure of the universe should be as perfect as it
goes. Correspondingly, the geometric space of quantum mechanics as suggested in
\cite{zil} emerges from a Zariski structure (see page~\pageref{zariskipage}) or
rather, from a sheaf of Zariski structures\footnote{The following three facts
clarify the connection between Zariski structures and the Isham-Butterfield
topos:
\begin{enumerate}
\item The sheaf of Zariski structures, the model of quantum mechanics, can be interpreted as a concrete realisation of an Isham-Butterfield topos.
\item The construction essentially generalises \cite{zil4} building a Zariski
  structure  corresponding to the non-commutative algebra represented by the
  canonical commutation relation $\mathrm{QP} -\mathrm{PQ}=i\hbar$. 
\item The analysis of the language and definability issues in the structure draws a clear line between notions which are {\em observable} (in the sense of physics) and which are not. 
\end{enumerate}
}.

It is equally important to note that the logical analysis inherent in our
method clarifies the correspondence between (possibly noncommutative) algebras
as they emerge in physics and geometry and the respective geometric spaces.  In
essence the algebras present us with the syntactic tools allowing to check in
calculations what can be seen graphically and dealt with geometrically. The
geometric space is thus a semantic interpretation of the syntactically given
data.  In classical cases, such as commutative finitely generated algebras,
this corresponds to the well-known duality at the foundation of algebraic
geometry. For  commutative  $C^*$-algebras we have the  Gel'fand-Naimark
duality linking those    to locally compact Hausdorff spaces. In
non-commutative cases the situation becomes much more complex but model theory
is in the best position to deal with the challenge. 

Another, different, line of collusion between categoricity and physics has been
explored by D.~Howard and I.~Toader in the past two decades
(see~\cite{How,Toad}). Their take on categoricity is more akin to the original
Veblen formulation than to the role categoricity has acquired in contemporary
model theory.

%On the other hand, if we ascribe to Galileo's motto that the physical world is 
%written in a mathematical language, it is natural to think that those theories in the core of mathematics, like logically perfect structures, should 
%play an important role in describing the ``reality''. \\

We finish this section with the conclusion that the principle of logical
perfection, as unconventional  as it may sound to some, does not disagree with
other modern approaches to the mathematical foundations of physics.

\section{Concluding remarks}

Our concept of logically perfect structures emerges as the result of a fifty
year classification project in logic. The theory is deep and technical but the
concept can be expressed in simple intuitive terms.

The defining property of logical perfection is uniqueness, or technically
uncountable categoricity. This property implies certain internal harmony:
homogeneity and the presence of a notion of dimension. This harmony is a
manifestation of a certain kind of geometricity, which itself is a consequence
of the infusion of geometric/topological ingredients in logic that brings forth
the flexibility and generality of logically perfect structures. Finally, since
logical structures are at the top of the classification hierarchy, they are
suitable as background structures for physics and represent a good idea of
geometric space in a very broad sense. 

An additional feature to support our notion of logically perfect structures is
the ``filtration'' of perfection provided by classification theory. As
mentioned above, classification theory not only places all first order theories
in a sort of map \emph{with respect to} categorical theories but provides a
kind of measure of going away from perfection. It provides technical ways to
measure, for arbitrary theories, what features of perfection might have been
lost and which ones remain. The second author's forthcoming interview with
Saharon Shelah explores further several peculiarities of this
connection\cite{ShInt}. \\ 

The features described above (uniqueness, geometricity, representability)
have concrete mathematical formulations, as we have briefly mentioned. In
addition, they help us to understand the role of those structures in the
wider program of studying the syntax/semantics duality. As we have tried to
show, logically perfect structures can be seen as located in the
geometric/semantical side of the mentioned duality, giving a new approach to
the notion of noncommutative (or quantum) geometric space, which traditionally
has been treated by means of syntactic/algebraic tools. Pursuing this program
of interpreting the duality between algebraic and geometric objects as an
extension of the 
duality between syntax and semantics appears to us as one of the most
interesting lines of research for the future, not only in mathematics. The idea
of representing one object by another (in this case its dual) can certainly be
extrapolated beyond mathematics. This idea deserves more investigation. \\

%Even though our motivation for formulating the idea of logical perfection is not
%philosophical, it is possible to see that there is a substantial philosophical
%component in that notion and that a study of logical perfection in more
%philosophical terms should give new insights that can be seen as an orthogonal
%complement to the technical approach. The development of the idea of Zariski
%structures and their generalizations has given birth to a ``new geometry'', a
%new ``representation theory'', and other new beautiful mathematical ideas that
%will lead to amazing discoveries in the future. This deserves to be studied
%by philosophers.

\end{document}